\documentclass{amsart}
\usepackage{amsmath,amsthm,amssymb,hyperref,amsfonts,epsf}

\newcommand{\al}{\alpha}
\newcommand{\de}{\delta}
\newcommand{\ep}{\epsilon}
\newcommand{\la}{\lambda}\newcommand{\La}{\Lambda}

\newcommand{\om}{\omega}
\newcommand{\Ga}{\Gamma}
\newcommand{\bR}{\mathbb{R}}\newcommand{\bZ}{\mathbb{Z}}
\newcommand{\bC}{\mathbb{C}}

\newcommand{\pt}{\partial_t}\newcommand{\pa}{\partial}

\newcommand{\D}{{\mathrm D}}\newcommand{\les}{{\lesssim}}

\newcommand{\beeq}{\begin{equation}}\newcommand{\eneq}{\end{equation}}

\newtheorem{thm}{Theorem}[section]
\newtheorem{prop}[thm]{Proposition}
\newtheorem{coro}[thm]{Corollary}
\newtheorem{rem}{Remark}[section]

\newtheorem{lem}{Lemma}

\newenvironment{prf}{\noindent {\bf Proof.} }{\endprf\par}
\def \endprf{\hfill  {\vrule height6pt width6pt depth0pt}\medskip}

\numberwithin{equation}{section}

\begin{document}

\title[Weighted Strichartz Estimates]{Weighted Strichartz Estimates with Angular Regularity and their
Applications}

\author{ Daoyuan Fang}
\address{Department of Mathematics, Zhejiang University, Hangzhou,
310027, China} \email{dyf@zju.edu.cn}

\thanks{The authors were supported by NSF of China 10571158 and 10871175.}

\author{Chengbo Wang}
\address{Department of Mathematics, Zhejiang University, Hangzhou,
310027, China}
\curraddr{Department of Mathematics, Johns Hopkins University, Baltimore,
Maryland 21218}
\email{wangcbo@jhu.edu}

\date{}
\dedicatory{} \commby{}

\begin{abstract}
In this paper, we establish an optimal dual version of trace estimate involving angular
regularity. Based on this estimate, we get the generalized Morawetz
estimates and weighted Strichartz estimates for the solutions to a
large class of evolution equations, including the wave and
Schr\"{o}dinger equation. As applications, we
prove the Strauss' conjecture with a kind of mild rough data for
$2\le n\le 4$, and a result of global well-posedness with small data
for some  nonlinear Schr\"{o}dinger equation with $L^2$-subcritical
nonlinearity.
\end{abstract}

\maketitle

\section{Introduction and Main Results}\label{6-sec-intro}
In the analysis of partial differential equations, e.g., the wave
and Schr\"{o}dinger equations, there are many results under the
assumption of the spherical symmetry. At the same time, there are
many radial estimates dealing with only the radial functions. In
general, such estimates or results would not hold for the general
case. A natural question is: how much additional angular regularity
of the functions is allowed to ensure that the radial results are
still valid in that case.

In principle, we believe that, for most of the results with radial
assumptions, there are the counterparts for the functions with
certain angular regularity. Recently, there have been some
interesting results in this direction, see e.g. Machihara-Nakamura-Nakanishi-Ozawa \cite{MaNaNaOz05},
Sterbenz \cite{Stbz05}, Kato-Nakamura-Ozawa \cite{KaNaOz07}, Cho-Ozawa \cite{Cho-Ozawa}. In this paper, we
make an attempt to get some systematic results in this direction.

To this end, we establish an optimal dual version of trace estimate involving angular
regularity. Based on this estimate, we get the generalized Morawetz
estimates and weighted Strichartz estimates for the solutions to a
large class of evolution equations, including the wave and
Schr\"{o}dinger equations. As the applications of our estimates, we
prove the Strauss' conjecture with a kind of mild rough data for
$2\le n\le 4$, and a result of global well-posedness with small data for
the nonlinear Schr\"{o}dinger equation.

Let $\mathcal{S}$ be the space of Schwartz function,
$\Delta_{\om}=\sum_{1\le i<j \le n} \Omega_{ij}^2$ be the
Laplace-Beltrami operator on $S^{n-1}\subset \bR^n$ with
$\Omega_{ij}=x_i
\partial_j-x_j \pa_i$, $\om\in S^{n-1}$, $\La_\om=\sqrt{1-\Delta_{\om}}$. For
$x\in \bR^n$, we introduce the polar coordinate $x=r\om$ with $r\ge
0$ and $\om\in S^{n-1}$. Let $\Delta=\sum_{i=1}^n
\pa_i^2=\pa_r^2+\frac{n-1}{r}\pa_r+\frac{1}{r^2}\Delta_\om$ be the
Laplacian and $\D=\sqrt{-\Delta}$. Based on the usual Besov spaces
$B^s_{p,q}$, we introduce the Besov spaces with angular regularity
as follows ($m\ge 0$)
$$B^{s,m}_{p,q,\om}=\La_\om^{-m}B^s_{p,q}=\{u\in B^s_{p,q} :
\| \La^m_\om u \|_{B^s_{p,q}}<\infty\}.$$ For the usual Sobolev
space $H^s=B^s_{2,2}$ or the homogeneous Besov space
$\dot{B}^s_{p,q}$, we can similarly define the spaces
$H^{s,m}_{\om}=B^{s,m}_{2,2,\om}$,
$\dot{B}^{s,m}_{p,q,\om}=\La_\om^{-m} \dot{B}^s_{p,q}$ and
$\dot{H}^{s,m}_{\om}=\dot{B}^{s,m}_{2,2,\om}$.

\vspace{0.3cm} \noindent$\bullet$ {\bf{Trace Lemma and Sobolev
Inequalities with Angular Regularity}}\vspace{0.3cm}

We first state the dual version of the trace lemma, which plays the central role in this
paper.
\begin{thm}[Dual Version of the Trace Lemma]\label{60-thm-trace}
Let $b\in (1,n)$ and $n\ge 2$, then we have the following equivalent
relation
  \beeq\label{60-est-trace}
        \|\La_\omega^{\frac{b-1}{2}} |x|^{-\frac{b}{2}}\widehat{g d \sigma}(x)\|_{L^2_x}
        \simeq
        \|g\|_{L^2_\omega}
  \eneq for any $g\in L^2_\omega$.
\end{thm}

Here $A \simeq B$ ($A \les B$) means that $c B\le A \le C B$ ($A \le C B$) with $C>c>0$, and, in what follows,
the constants $C$ and $c$ might change at each occurrence.

From \eqref{60-est-trace}, we know that
  \beeq\label{60-est-trace2}
        \||x|^{-\frac{b}{2}}\widehat{g d \sigma}(x)\|_{L^2_x}
        \les
        \|g\|_{H^{s}_\omega}
  \eneq if  $s\ge \frac{1-b}{2}$ and $b\in(1,n)$.

\begin{rem}
The estimate \eqref{60-est-trace2} has been proved by Wang in
\cite{Wang91} with $s=0$ and $b\in(1,n)$ , and Hoshiro
\cite{Hoshiro97} with $s=\frac{1-b}{2}$ and $n\ge 3$. The novelty of
our estimate \eqref{60-est-trace} consists in getting the equivalent
relation instead of the usual inequality.\end{rem}

\begin{rem}
For the estimate \eqref{60-est-trace} or \eqref{60-est-trace2}, the
condition on $b$ is also necessary, see Remark \ref{6-rem-3.1}.
\end{rem}

By duality, scaling and the Sobolev embedding in $S^{n-1}$
$$H^{s}_\om\subset L^\infty_\om\textrm{ with }s>\frac{n-1}{2},$$
we can get the following Sobolev estimates which are of independent
interest.
\begin{coro}\label{60-thm-Sobolev-trace}
  Let $b\in(1,n)$ and $n\ge 2$, we have
  \beeq\label{60-est-Sobolev-infty} \sup_{r>0}\ r^{\frac{n-b}{2}} \|f(r\omega)\|_{L^2_\omega}
  \les \| \D^{\frac{b}{2}} \La_\omega^{\frac{1-b}{2}}
  f\|_{L^2_x} \eneq and
  \beeq\label{60-est-Sobolev-dual}
   \| \D^{-\frac{b}{2}} \La_\omega^{\frac{b-1}{2}}
  f\|_{L^2_x}\les
  \| |x|^{\frac{b-n}{2}} f(x)\|_{
L^1_{|x|^{n-1} d |x|} L^2_\omega} \eneq for any $f\in \mathcal{S}$.
Moreover, if $s>\frac{n-b}{2}$, then
\beeq\label{60-est-Sobolev-infty2} \sup_{r>0}\  r^{\frac{n-b}{2}}
\|f(r\omega)\|_{L^\infty_\omega}
  \les \| \D^{\frac{b}{2}} \La_\omega^{s}
  f\|_{L^2_x}.
\eneq
\end{coro}

\begin{rem}
When $f$ is a radial function, the estimate
\eqref{60-est-Sobolev-infty} has been proved essentially by Li-Zhou
(Theorem 2.10 in \cite{LiZhou95}). The estimate
\eqref{60-est-Sobolev-infty2} for radial $f$ with $b = 2$ and $n \ge
3$ reduces to Ni's inequality \cite{Ni82}.
\end{rem}
\begin{rem}
Recently, Cho-Ozawa prove the estimate
\eqref{60-est-Sobolev-infty2} with $s>n-1-\frac{b}{2}$ in
\cite{Cho-Ozawa}. Our result improves the angular regularity.
Moreover, the requirement of angular regularity in our estimate is
essentially optimal, since the index
$\frac{n-b}{2}+\frac{b}{2}=\frac{n}{2}$ is precisely the infimum of $s$
so that $H^s\subset L^\infty$.
\end{rem}

As a side remark, we have the following alternative trace estimate
for the forbidden case $b=1$.
\begin{prop}
\label{60-thm-trace-v2} For any compact $C^\infty$ hypersurface
$M\subset \bR^n$ and $g \in L^2_M$, we have
  \beeq\label{60-est-trace-v2}
        \sup_{x_0, R} R^{-\frac{1}{2}} \|\widehat{g d M}(x)\|_{L^2_{B(x_0,R)}}
        \les
        \|g\|_{L^2_M}
        \les
        \limsup_{x_0, R\rightarrow \infty} R^{-\frac{1}{2}} \|\widehat{g d M}(x)\|_{L^2_{B(x_0,R)}}.
  \eneq
Moreover, by duality, we have $\dot{B}^{\frac{1}{2}}_{2,1}\subset
L_M^2 $, and by rescaling, we get
  \beeq\label{60-est-Sobolev-infty-v2} r^{\frac{n-1}{2}} \|f(r y)\|_{L^2_{y\in M}}
  \les \|   f\|_{\dot{B}^{\frac{1}{2}}_{2,1}} \eneq  for any $f\in \mathcal{S}$.
\end{prop}
\begin{rem}
The inequality \eqref{60-est-trace-v2} is precisely Theorem 2.1 and Theorem 2.2 of
Agmon-H\"{o}rmander \cite{AgmHor76}.
\end{rem}

Let $\Delta_M$ be the Laplace-Beltrami operator on $M$ and
$\La_M=(1-\Delta_M)^{\frac{1}{2}}$, then, from
\eqref{60-est-Sobolev-infty-v2}, we can get the following estimate,
\beeq\label{60-est-Sobolev-infty2-v2} r^{\frac{n-1}{2}} \|f(r
y)\|_{L^\infty_{y\in M}}
  \les \|  \La_M^{s} f\|_{\dot{B}^{\frac{1}{2}}_{2,1}}\eneq with
  $s>\frac{n-1}{2}$.
This also relaxes the condition $s>n-\frac{3}{2}$ of Cho-Ozawa
\cite{Cho-Ozawa} in the case of $M=S^{n-1}$.

As an application of \eqref{60-est-Sobolev-infty2} and
\eqref{60-est-Sobolev-infty2-v2}, we can get the following result of
compact embedding (see e.g. p7 of \cite{Cho-Ozawa} for the proof,
and see e.g. Section 1.7 of Cazenave \cite{Cazenave} for the
previous radial result)
\begin{prop}[Compact Embedding]\label{6-thm-Sobolev-Cpt}
The embedding $H^{\frac{b}{2}, m}_\om \subset L^p$ is compact for
$b\in (1,n)$, $m>\frac{n-b}{2}$ and $2<p<\frac{2 n}{n-b}$. Moreover,
the embedding $B^{\frac{1}{2},m}_{2,1,\om} \subset L^p$ is compact
for $m>\frac{n-1}{2}$ and $2<p<\frac{2 n}{n-1}$.
\end{prop}

\vspace{0.3cm} \noindent$\bullet$ {\bf{Generalized Morawetz
Estimates}}\vspace{0.3cm}

It is well-known that we can get certain generalized Morawetz
inequality (or the local smoothing effect for the dispersive
case $a>1$) from the knowledge of the inequality like \eqref{60-est-trace2}.
Here we can give a more refined estimate because of the improved version of the trace lemma.

\begin{thm}\label{60-thm-Morawetz} If $b\in(1,n)$ and $a>0$, then we
have \beeq\label{60-est-Morawetz}\||x|^{-\frac{b}{2}} e^{i t \D^a}
f\|_{L^2_{t,x} }\simeq
\|\D^{\frac{b-a}{2}}\La_\omega^{\frac{1-b}{2}}f\|_{L^2_x}\eneq  for
any $f\in \mathcal{S}$. Moreover, for the endpoint case $b=1$, we
have the following local estimate \beeq\label{60-est-Morawetz-Local}
\sup_{x_0, R} R^{-\frac{1}{2}} \| e^{i t \D^a} f\|_{L^2_{t, B(x_0,
R) } } \les \|\D^{\frac{1-a}{2}} f\|_{L^2_x} \eneq  for any $f\in
\mathcal{S}$.
\end{thm}

The generalized Morawetz estimates usually take the form \beeq
\label{60-est-Morawetz2} \| |x|^{-\frac{b}{2}} e^{i t \D^a} f
\|_{L^2_{t,x} } \les \|\D^{\frac{b-a}{2}} \La_\omega^{s} f\|_{L^2_x}
. \eneq The novelty of our estimate \eqref{60-est-Morawetz} consists
in getting the equivalent relation instead of the usual inequality.
Note that a similar phenomenon was found recently by Vega-Visciglia
\cite{VeVi07} to the estimate \eqref{60-est-Morawetz-Local} for the
Schr\"{o}dinger equation ($a=2$).

\begin{rem} The estimates \eqref{60-est-Morawetz2} without angular smoothing
index (s=0) have been proved by many authors for the special case
$a=1,2$. Morawetz \cite{Mor68} first got the estimate with $s=0$,
$b=3$ and $n\ge 4$ for the wave equation ($a=1$). Kato and Yajima
\cite{KatoY88} gave the estimate with $s=0$, $b\in (1,2]$ and $n\ge
3$ for the Schr\"{o}dinger equation ($a=2$). Sugimoto \cite{Su98}
and Vilela \cite{Vilela01} prove the estimate for $a=2$ with $s=0$.
For the estimate with an angular smoothing index, Hoshiro
\cite{Hoshiro97} and Sugimoto \cite{Su98} get some of the estimates
for $a=1, 2$ with $s=\frac{1-b}{2}$.
\end{rem}

\begin{rem}
The local smoothing estimate \eqref{60-est-Morawetz-Local} was
obtained by Kenig-Ponce-Vega for $a=2$ in \cite{KPV91}. As far as we
know, the space-localized estimate with $a=1$ was first obtained by
Smith-Sogge in Lemma 2.2 of \cite{SmSo00}.\end{rem}

\begin{rem}
In the case of wave equation ($a=1$), there also exists an important
time-localized estimate, which was first obtained by
Keel-Smith-Sogge \cite{KSS02}
$$\log(2+T)^{-\frac{1}{2}}\| \langle x\rangle^{-\frac{1}{2}} e^{i t \D}
f\|_{L_{[0,T]}^2 L^2_{x}}\les \|f\|_{L^2_x}.
$$
\end{rem}

\vspace{0.3cm} \noindent$\bullet$ {\bf{Weighted Strichartz
Estimates}}\vspace{0.3cm}

For the operator $e^{i t \D^a}$ with $a>0$, the so-called Strichartz
estimates usually take the form (see e.g. Keel-Tao \cite{KeTa98})
\beeq\label{60-est-Stri-Gene}\|e^{i t \D^a} f\|_{ L^q_t L^r_{x} }
\les \| f\|_{\dot{H}^s_x},\eneq  where
$s=\frac{n}{2}-\frac{a}{q}-\frac{n}{r}$ by scaling. The Strichartz
estimates have been proved to be very useful in the study of the
well-posed problems, see e.g. \cite{Cazenave}, \cite{LdSo95}, \cite{ShSt00} and \cite{Sog95}.
In practice, it is
interesting and meaningful to consider the generalized Strichartz
estimates of the following type \beeq\label{60-est-se1} \|e^{i t
\D^a} f\|_{L^q_t L^r_{|x|^{n-1}d|x|} L^p_\omega }\les \|
f\|_{\dot{H}^{s,s_1}_{\om}}.\eneq Moreover, we are interested in
obtaining the weighted Strichartz estimates
$$\||x|^{-\al} e^{i t
\D^a} f\|_{L^q_t L^r_{|x|^{n-1}d|x|} L^p_\omega }\les \|
f\|_{\dot{H}^{s,s_1}_{\om}}.$$

If we interpolate between the estimates \eqref{60-est-Morawetz} and
the Sobolev inequality \eqref{60-est-Sobolev-infty}, we can get the
following weighted Strichartz estimates.
\begin{thm}[Weighted Strichartz Estimates]\label{60-thm-Strichartz-weighted}
If $b\in (1,n) $, $a>0$ and $r\in [2,\infty]$, we have
\beeq\label{60-est-Morawetz-Strichartz}
\||x|^{\frac{n}{2}-\frac{n}{r}-\frac{b}{2}} e^{i t \D^a} f (x)
\|_{L^r_{t, |x|^{n-1}d |x| } L^2_\omega } \les
\|\D^{\frac{b}{2}-\frac{a}{r}}\La_\omega^{\frac{1-b}{2}}f\|_{L^2_x},\eneq
for any $f\in \mathcal{S}$.
\end{thm}
The estimates stated in Theorem \ref{60-thm-Strichartz-weighted} is
the homogeneous estimates. In practice, it is often important to
give the inhomogeneous estimates. By the Christ-Kiselev
lemma (Theorem 1.2 in \cite{ChKi01}), we can get the inhomogeneous
estimates. In conclusion, we have
\begin{thm}\label{60-thm-StriInhom}
Let $q, \tilde{q}\in [2, \infty]$, $\frac{n}{q}-\al,
\frac{n}{\tilde{q}}-\tilde{\al}\in (0, \frac{n-1}{2})$,
$s=\frac{n+a}{q}-\frac{n}{2}-\al$,
$s_1=\frac{n-1}{2}+\al-\frac{n}{q}$ (note that
$s+s_1=\frac{a}{q}-\frac{1}{2}$), and $\tilde{s}, \tilde{s_1}$
similarly defined. Then we have
\beeq\label{60-est-Morawetz-Strichartz2} \||x|^{-\al} \D^s
\La_\om^{s_1} e^{i t \D^a} f (x) \|_{L^q_{t, |x|^{n-1}d |x| }
L^2_\omega } \les \|f\|_{L^2_x},\eneq and by duality, one can get
\beeq\label{60-est-Mora-Stri-dual} \| \int \D^s \La_\om^{s_1} e^{-i
s \D^a} F (s,x) d s \|_{L^2_x} \les \| |x|^{\al} F\|_{L^{q'}_{t,
|x|^{n-1}d |x| } L^2_\omega }.\eneq Moreover, we have the following
inhomogeneous estimates \beeq\label{60-est-Mora-Stri-dual2}\|
\int_0^t \D^s \La_\om^{s_1} e^{i (t-s) \D^a} F (s,x) d s
\|_{L^\infty_t L^2_x} \les \| |x|^{\al} F\|_{L^{q'}_{t, |x|^{n-1}d
|x| } L^2_\omega },\eneq and \beeq\label{60-est-Mora-Stri-Inhom} \|
\int^t_0 |x|^{-\al} \D^{s+\tilde{s}} \La_\om^{s_1+\tilde{s_1}} e^{i
(t-s) \D^a} F (s,x) d s \|_{L^q_{t, |x|^{n-1}d |x| } L^2_\omega }
\les \| |x|^{\tilde{\al}} F\|_{L^{\tilde{q}'}_{t, |x|^{n-1}d |x| }
L^2_\omega }\eneq with $q>\tilde{q}'$.
\end{thm}

\begin{rem}
Harmse \cite{Harmse} and Oberlin \cite{Ober89} got the inhomogeneous
inequality \beeq\label{60-est-Harmse} \|w\|_{L^q_{t,x}}\les
\|F\|_{L^r_{t,x}}, \eneq for the solution $w$ of the equation
$(\pt^2-\Delta) w=F$ with data $(0,0)$, if
$\frac{n+1}{r}-\frac{n+1}{q}=2$ and
\beeq\label{60-est-Harmse-cond}\frac{n+1}{2n }
-\frac{2}{n+1}<\frac{1}{q}<\frac{n-1}{2n}.\eneq Our estimate
\eqref{60-est-Mora-Stri-Inhom} generalizes the Harmse-Oberlin
estimate for $a=1$ to the general cases with weights.
\end{rem}

In particular, if we choose $b\in (1,n)$ such that
$\frac{n}{2}-\frac{n}{r}-\frac{b}{2}=0$ in Theorem
\ref{60-thm-Strichartz-weighted}, we can get the following
generalized Strichartz estimates with $q=r$ in presence of angular
regularity.
\begin{coro}\label{60-thm-Strichartz}
Let $a>0$, $r\in (\frac{2 n}{n-1},\infty)$ and $p\in [2,\infty)$, we
have \beeq\label{60-est-Strichartz-Angu-q--r} \| e^{i t \D^a} f (x)
\|_{L^r_{t, |x|^{n-1}d |x| } L^p_\omega } \les
\|\D^{\frac{n}{2}-\frac{n+a}{r}}\La_\omega^{\frac{n}{r}-\frac{n-1}{p}}
f\|_{L^2_x},\eneq for any $f\in \mathcal{S}$.
\end{coro}

\vspace{0.3cm} {\noindent $\bullet$ \bf Generalized Strichartz
Estimates for the Wave Equation}\vspace{0.3cm}

In the case of the wave equation ($a=1$), it is well known that we
have the classical Strichartz estimates \eqref{60-est-Stri-Gene}
(see \cite{FW2}) if \beeq\label{60-est-wave-Stri}\frac{1}{q}\le
\min(\frac{1}{2},\frac{n-1}{2}(\frac{1}{2}-\frac{1}{r})), (q,r)\neq
(\max(2,\frac{4}{n-1}),\infty), (q,r)\neq (\infty,\infty). \eneq The
result in Corollary \ref{60-thm-Strichartz} extends the Strichartz
estimates to the case of $q=r<\infty$ and
$$\frac{1}{q}<(n-1)(\frac{1}{2}-\frac{1}{r}).$$ Thus it is natural
to guess that there is a similar result in the more general case
$q\neq r$. It is in fact the case at least for $r = p$ (see Sterbenz
\cite{Stbz05} for $n\ge 4$, and Section \ref{6-sec-Strichartz} for the full range $n\ge 2$).
\begin{thm}\label{60-thm-Strichartz-Sterbenz}Let $s=n(\frac{1}{2}-\frac{1}{r})-\frac{1}{q}$,
\beeq\label{60-est-Conj1-wave-skn}s_{kn}=\frac{2}{q}-(n-1)(\frac{1}{2}-\frac{1}{r}),\eneq
and
\beeq\label{60-est-Conj1-wave}\frac{n-1}{2}(\frac{1}{2}-\frac{1}{r})<\frac{1}{q}<(n-1)(\frac{1}{2}-\frac{1}{r}),\
q\ge 2,\eneq then we have the estimates
$$\|e^{i t \D} f\|_{L^q_t L^r_{x}}\les \|
f\|_{\dot{H}^{s,s_1}_{\om}}$$ for any $s_1>s_{kn}$.
\end{thm}

In our previous paper \cite{FW2}, we prove that the endpoint
Strichartz estimates
$$\|e^{i t \D}f\|_{L^q_t L^r_x}\les \|f\|_{\dot{H}^s}$$
with $(q,r)=(4,\infty)$ and $n=2$ does not hold in general, and
holds for radial functions. As a corollary of
Theorem~\ref{60-thm-Strichartz-Sterbenz}, we can recover this
estimate by interpolating the angular $L^{4-}L^\infty$ estimate and
the classical $L^{4+}L^\infty$ estimate, if we add certain angular
regularity. Thus, by  combining it with the $n=3$ result of
\cite{MaNaNaOz05}, we have the following result.
\begin{coro}
Let $n=2,3 $, for any $\ep>0$, we have
\beeq\label{60-est-Enpt-Stri-Ang} \|e^{i t\D}
f(x)\|_{L_t^{\frac{4}{n-1}} L_x^\infty}\les  \|f\|_{\dot
H^{\frac{n+1}{4},\ep}_{\om}}.\eneq\end{coro}

\vspace{0.3cm} {\noindent $\bullet$ \bf Generalized Strichartz
Estimates for the Schr\"{o}dinger Equation}\vspace{0.3cm}

 The case $a=2$ is just the case of
the Schr\"{o}dinger equation. In this case, recall that we have the
classical Strichartz estimates (see e.g. Keel-Tao \cite{KeTa98})
$$\|e^{i t
\Delta} f\|_{L^q_t L^r_x}\les \|f\|_{L^2},$$ if $$
\frac{1}{q}=\frac{n}{2}(\frac{1}{2}-\frac{1}{r})\le \frac{1}{2},
(q,r,n)\neq (2,\infty,2).$$ As in the case of the wave equation (see
\cite{FW2}), we can generalize the estimates to
\beeq\label{60-est-Stri-Schro}\|e^{i t \Delta}f\|_{L^q_t L^r_x}\les
\|\D^{\frac{n}{2}-\frac{2}{q}-\frac{n}{r}} f\|_{L^2_x},\eneq for
\beeq\label{60-est-Schro}\frac{1}{q}\le
\min(\frac{1}{2},\frac{n}{2}(\frac{1}{2}-\frac{1}{r})), (q,r)\neq
(\infty,\infty), (2,\infty).\eneq

Then, the estimates \eqref{60-est-Stri-Schro} require $\frac{2
(n+2)}{n}\le r<\infty$ in the case of $q=r$. However, by adding some
additional angular regularities, we can relax this restriction to
$r>\frac{2 n}{n-1}$ in \eqref{60-est-Strichartz-Angu-q--r}.
Moreover, by interpolate with the known Strichartz estimates for
$q=r=\frac{2 (n+2)}{n}$, we can improve the estimate
\eqref{60-est-Strichartz-Angu-q--r} with $a=2$ further.
\begin{coro}\label{60-thm-Strichartz-Schro} Let $n>2$ and
$$r\in (\frac{2 n}{n-1}, \frac{2 (n+2)}{n}),$$
we have \beeq\label{60-est-Stri-Schro-Angu}\|e^{i t \Delta} f
\|_{L^r_{t, x}}\les \|\D^{\frac{n}{2}-\frac{n+2}{r}}
\La_{\om}^{\frac{n-1}{n-2}(\frac{n+2}{r}-\frac{n}{2})+\ep}
f\|_{L^2_x},\eneq for any $\ep>0$ and $f\in \mathcal{S}$.
\end{coro}

\begin{rem}
Recently, when the data $f$ is radial, Shao \cite{Shao} generalizes
the estimates \eqref{60-est-Stri-Schro} to the range of $2
\frac{2n+1}{2n-1}<q=r<2 \frac{n+2}{n}$.
\end{rem}

In general, we {\bf conjecture} that, the estimate
\eqref{60-est-se1} with $a=2$ and $s=\frac{n}{2}-\frac{2}{q}-\frac{
n}{r}$ holds true with $s_1>s_{kn}$ for
\beeq\label{60-est-Conj1-Schro}
\frac{n}{2}(\frac{1}{2}-\frac{1}{r})<\frac{1}{q} <\frac{2
n-1}{2}(\frac{1}{2}-\frac{1}{r}) \textrm{ and } q\ge 2,\eneq where
\beeq\label{60-est-Conj1-Schro-skn}s_{kn}=\frac{2}{q}+\frac{2
n-1}{r}-\frac{n-1}{p}-\frac{n}{2}.\eneq It may be interesting to
note that in the case of $p=r$, $s_{kn}=\frac{2}{q}+\frac{
n}{r}-\frac{n}{2}$, which is just $-s$.

Next, for the above estimates, we give two applications.

\vspace{0.3cm} \noindent$\bullet$ {\bf{Applications to the Wave
Equation}}\vspace{0.3cm}

Let $x\in\bR^n$ with $n\ge 2$, $F_p(u)=\la |u|^p$ ($\la\in \bR
\backslash \{0\}$, $p>1$), $s_c=\frac{n}{2}-\frac{2}{p-1}$,
$s_{sb}=\frac{1}{2}-\frac{1}{p}$, $p_{conf}=1+\frac{4}{n-1}$,
$p_h=1+\frac{4 n}{(n+1)(n-1)}$  and $p_c$ be the solution of the
quadratic equation
$$ (n-1) p_c^2 - (n+1) p_c - 2 = 0,\ p_c>1. $$ Note that $p>p_c$ if and only if $s_c>s_{sb}$. Consider the
following semi-linear wave equations for $u:
\bR\times\bR^n\rightarrow \bR$, \beeq \label{60-eqn-SLW}
\left\{\begin{array}{l} (\pt^2-\Delta ) u = F_p (u)\\
u(0,x)=f, \pt u(0,x)=g.
\end{array}\right. \eneq
Strauss' conjecture asserts that the problem \eqref{60-eqn-SLW} has
a global solution for $p>p_c$, when the initial data $(f,g)$ is
sufficiently small and smooth with compact support. This conjecture
was finally completed by Georgiev-Lindblad-Sogge in
\cite{GLS97} (see also Tataru \cite{Ta01} for another proof). In \cite{GLS97}, the authors
raised an interesting problem: ``{\sf Under what kind of the low
regularity assumptions on the data, the Strauss conjecture still
holds true?}"

When $p\ge p_{conf}$, Lindblad-Sogge \cite{LdSo95} have
succeeded in getting the global well-posed result for
\eqref{60-eqn-SLW} with small data $(f,g)\in \dot{H}^{s_c}\times
\dot{H}^{s_c-1}$, where $s_c$ is the minimal regularity assumption.
When the data $(f,g)\in \dot{H}^{s_c}\times \dot{H}^{s_c-1}$ are
small and radial, there are results dealing with either $n\le 4$ or
$p>p_h$, see e.g. Sogge \cite{Sog95}, Lindblad-Sogge \cite{LdSo95}
and Hidano \cite{Hidano}.

We can apply  the above estimates to this problem for the small data
$(f,g)\in \dot{H}^{s_c,s_1}_\om \times \dot{H}^{s_c-1, s_1}_\om$
with some $s_1$. Following the arguments in Section 8 of
Lindblad-Sogge \cite{LdSo95}, we will first prove the following
\begin{thm}\label{60-thm-LdSo95}
Let $n\ge 2$, $p_h<p<p_{conf}$ (i.e.,
$\frac{1}{2n}<s_c<\frac{1}{2}$) and $s_1>\frac{1}{2}-s_c$. Suppose
that $$(f,g)\in \dot{H}^{s_c,s_1}_\om \times \dot{H}^{s_c-1,
s_1}_\om$$ with small enough norm, then there is a unique global
weak solution $u$ to \eqref{60-eqn-SLW} satisfying
$$u\in C_t \dot{H}_x^{s_c}\cap C^1_t
\dot{H}_x^{s_c}\cap L^{q}_{t,x} \textrm{ with }
q=\frac{(n+1)(p-1)}{2}.$$
\end{thm}

Moreover, we can prove the Strauss' Conjecture with a kind of mild
rough data for $n\le 4$ in the sense of the following theorem.
\begin{thm}\label{60-thm-Hidano}
Let $2\le n\le 4$, $p_c<p<p_{conf}$ and $s_1= \frac{1}{p-1}$.
Suppose that
$$(f,g)\in \dot{H}^{s_c,s_1}_\om \times \dot{H}^{s_c-1, s_1}_\om$$ with small enough norm,
then there is a unique global weak solution $u\in C_t
\dot{H}_{\om}^{s_c,s_1}\cap C^1_t \dot{H}_{\om}^{s_c-1,s_1}$ to
\eqref{60-eqn-SLW} satisfying
$$|x|^{-\al} u\in L^{p}_{t,|x|^{n-1}d|x|}
H^{s_2}_\om,$$ for $\al=\frac{n+1}{p}-\frac{2}{p-1}$ and
$s_2=s_1+s_c-s_{sb}$.
\end{thm}

\begin{rem}
At the final stage of preparation, we learned from Professor Sogge that they have independently
obtained a related result for the equation on the exterior domain (see Hidano-Metcalfe-Smith-Sogge-Zhou \cite{HMSSZ}).
\end{rem}

\vspace{0.3cm} \noindent$\bullet$ {\bf{Applications to the
Schr\"{o}dinger Equation}}\vspace{0.3cm}

Let $x\in\bR^n$ with $n\ge 1$, $F_p(u)=\la |u|^p$ or  $F_p(u)=\la
|u|^{p-1}u$ ($\la\in \bR \backslash \{0\}$, $p>1$),
$s_c=\frac{n}{2}-\frac{2}{p-1}$, $p_{L2}=1+\frac{4}{n}$,
$p_l=1+\sqrt{\frac{2}{n-1}}$. Consider the following nonlinear
Schr\"{o}dinger equation for $u: \bR\times\bR^n\rightarrow \bC$,
\beeq \label{60-eqn-SLS}
\left\{\begin{array}{l} (i \pt-\Delta ) u = F_p (u)\\
u(0,x)=f.
\end{array}\right. \eneq

Cazenave and Weissler \cite{CaWe90} proved that if $s_c\ge 0$ (i.e.
$p\ge p_{L2}$) and $[s_c]<p-1$ ($s_c<p$ if $s_c\in \bZ$), this
problem is global well-posed in $C_t H_x^{s_c}$ (for $f\in
H_x^{s_c}$ with small enough $\dot{H}_x^{s_c}$ norm). Moreover, if
$s_c<0$, the problem is global well-posed in $C_t L_x^2$, and fail
to be uniformly well-posed in $C_t H_x^s$ for any $s<0$ (see
Birnir-Kenig-Ponce-Svanstedt-Vega \cite{BKPSV} or Christ-Colliander-Tao \cite{ChCoTa03}). Recently, by
assuming $f\in \dot{H}^{s_c}$ to be small and radial, Hidano
\cite{Hidano07} get a global result for some $L^2$-subcritical
nonlinearity $1+\frac{4}{n+1}<p<p_{L2}$ and $n\ge 3$. We generalize
his result to the angular case as follows.

\begin{thm}\label{60-thm-SLS}
Let $3\le n\le 6$, $p_l<p<p_{L2}$ and $s_1= \frac{1}{p-1}$. Suppose
that $f\in \dot{H}^{s_c,s_1}_\om$ with small enough norm, then there
is a unique global weak solution $u\in C_t \dot{H}^{s_c,s_1}_\om$ to
the equation \eqref{60-eqn-SLS}, satisfying
$$|x|^{-\al} u\in L^{q}_{t,|x|^{n-1}d|x|}
H^{s_2}_\om$$ for $\al=\frac{n+2}{q}-\frac{2}{p-1}$ and
$s_2=\frac{n-1}{2} + \frac{2}{q} - \frac{1}{p-1}$, where $q$ satisfy
the restriction
$$\frac{2}{q}\in [\frac{1}{p}, 1]\cap (\frac{2}{p-1}-\frac{n-1}{2},
\frac{2}{p-1}-\frac{n+1}{2 p}) \cap [\frac{1}{p-1}-\frac{n-1}{2 p},
\frac{1}{p-1}-\frac{n-3}{2 p}).$$
\end{thm}

This paper will be organized as follows. In Section
\ref{6-sec-basic}, we collect some preliminary results concerning on
the analysis on the sphere and the knowledge of the Bessel
functions. In Section \ref{6-sec-traceStri}, we prove the dual version of the trace
lemma (Section \ref{6-sec-trace}), the generalized Morawetz
estimates (Section \ref{6-sec-Morawetz}) and the generalized
Strichartz estimates (Section \ref{6-sec-Strichartz}), in the
presence of angular regularity. In Section \ref{6-sec-SLW}, we give
the first application of the estimates obtained in Section
\ref{6-sec-traceStri} to the Strauss' conjecture with a kind of mild
rough data for $n\le 4$ or $p_h<p<p_{conf}$. In the final Section
\ref{6-sec-SLS}, we give another application to a result of global
well-posedness with small data for the nonlinear Schr\"{o}dinger
equation.

\section{Preliminary}\label{6-sec-basic}

For any $x\in \bR^n$, we introduce the polar coordinates $r=|x|$ and
$\om\in S^{n-1}$ such that $x=r\om$. Let $\langle x
\rangle=\sqrt{1+|x|^2}$ and $H(t)$ be the usual Heaviside function
($H(t)=1$ if $t\ge 0$ and $H(t)=0$ else). For a set $E$, we use
$|E|$ to stand for the measure or cardinality of the set $E$
depending on the context.

The proof of the trace lemma is based on the expansion of a function
defined on the sphere with respect to the spherical harmonics. Here
we describe the expansion precisely.

Let $n\ge 2$. For any $k\ge 0$, we denote by $\mathcal{H}_k$ the
space of spherical harmonics of degree k on $S^{n-1}$, by
$d(k)=\frac{2k+n-2}{k}C^{n+k-3}_{k-1}\simeq \langle k\rangle^{n-2}$
its dimension, and by $\{Y_{ k, 1} , \cdots , Y_{ k, d(k)} \}$ the
orthonormal basis of $\mathcal{H}_k$. It is well known that
$L^2(S^{n-1}) = \bigoplus_{k=0}^\infty \mathcal{H}_k$ and that $F(t,
x) = F(t, r \om)$ has the expansion
\beeq\label{6-est-Expansion-SHarm}F(t, r\om)=\sum_{k=0}^\infty
\sum_{l=1}^{d(k)}  a_{k,l}(t,r) Y_{k,l}(\omega).\eneq By
orthogonality, we observe that $\|F(t,
r\cdot)\|_{L^2_\om}=\|a_{k,l}(t,r)\|_{l^2_{k,l}}$.

Let $\Delta_{\om}$ be the Laplace-Beltrami operator on $S^{n-1}$
and $\La_\om=\sqrt{1-\Delta_{\om}}$. Then we have $\Delta_\om
Y_{k,l}=-k(k+n-2)Y_{k,l}$. Based on this fact, we naturally
introduce the Sobolev space $H^s_\om=\La_\om^{-s}L^2_\om$ on the
sphere $S^{n-1}$ and we have
\beeq\label{6-est-Sobolev-Sphere}\|F(t,r\cdot)\|_{H^s_\om}=\|\La_\om^s
F(t,r\om)\|_{L^2_\om}\simeq \|\langle k\rangle^s a_{k,l}(t,r)
\|_{l^2_{k,l}}.\eneq

The nonlinear estimates on the sphere, such as the Sobolev
embedding, the Leibniz rule and the Moser estimate, easily transfers
from the Euclidean case (see e.g. \cite{KaNaOz07}). In particular,
we have the following Moser estimate (for the Euclidean case, c.f.
Kato \cite{Kato95}) \beeq\label{6-est-Moser-Sphere}\|\La_\om^s
F_k(u)\|_{L^r_\om}\les \|u\|_{L^q_\om}^{k-1} \|\La_\om^s
u\|_{L^p_\om} \eneq for $s\in [0,m]$ and $p,q,r\in(1,\infty)$ with
$\frac{1}{r}=\frac{k-1}{q}+\frac{1}{p}$, where $F_k\in C^m$ with
$$F(0)=0, |\partial^\al F(x)|\les |x|^{k-|\al|}, 1\le |\al|\le m\le
k.$$

The spherical harmonics are closely connected with the special
functions such as the Gamma function, Bessel functions and so on.
Let
$$\Ga(s)=\int^\infty_0 e^{-r} r^{s-1} d r, s>0$$ be the Gamma function,
the Bessel function of order $k>-\frac{1}{2}$ is defined by
\beeq\label{6-est-Bessel-Def} J_k(t) = \frac{2^{-k} t^k }{
\Ga(k+\frac{1}{2}) \Ga(\frac{1}{2})} \int_{-1}^1 e^{i t s}
(1-s^2)^{k-\frac{1}{2}} d s. \eneq

For these functions, we have some well-known asymptotics. For the
Gamma function, we have the Stirling's formula (\cite{G}, p.421)
\beeq\label{6-est-Gamma-asymptotic}\Ga(t)\simeq \sqrt{2 \pi}
t^{t-\frac{1}{2}}e^{-t},\textrm{ as }t\rightarrow \infty.\eneq For
the Bessel function, we have
\beeq\label{6-est-Bessel-asymptotic}J_k(t)\simeq
\sqrt{\frac{2}{\pi}} t^{-\frac{1}{2}}\cos(t-\frac{2 k+1}{4}\pi)
\textrm{ as } t\rightarrow \infty, J_k(t)\simeq 2^{-k} t^{k}
\textrm{ as } t\rightarrow 0.\eneq

The following result is proved to be very useful. (Note that in
their notation, $P(x)=|x|^k Y_{k,l}(\omega)$)
\begin{lem}[IV. Theorem 3.10 in Stein-Weiss \cite{SW71}] Let
\beeq\label{6-est-gkl-from-fkl}
  g_{k,l}(\rho)=(2 \pi)^{\frac{n}{2}} i^{- k} \int^\infty_0 f_{k,l}(r)
  \rho^{-\frac{n-2}{2}}
        J_{k+\frac{n-2}{2}}(r \rho) r^{\frac{n}{2}} d r,
\eneq then we have \beeq\label{6-est-F-SphericalHarmonic}
  \widehat{f_{k,l}(r)Y_{k,l}(\omega)}(\xi) =
  g_{k,l}(|\xi|)Y_{k,l}(\frac{\xi}{|\xi|})
\eneq
\end{lem}

As a corollary, if we set $f(r)=\de(r-1)$, then we get that
\beeq\label{6-est-F-SphericalHarmonic2}\widehat{Y_{k,l}(\omega)
d\sigma(\omega)}(\xi) =(2 \pi)^{\frac{n}{2}} i^{- k}
  |\xi|^{-\frac{n-2}{2}}
        J_{k+\frac{n-2}{2}}(|\xi|)Y_{k,l}(\frac{\xi}{|\xi|}).\eneq
In particular, if $k=0$, then we have the well known
\beeq\label{6-est-F-d-sigma}\widehat{d \sigma}(\xi) =(2
\pi)^{\frac{n}{2}}  |\xi|^{-\frac{n-2}{2}}
        J_{\frac{n-2}{2}}(|\xi|).\eneq

We will also need the following special form of the
Weber-Schafheitlin integral formula for Bessel functions (p.403 of
Watson \cite{Watson44}).
\begin{lem}
Let $\mu, \nu, \la\in \bR$ such that $\mu+\nu+1>\la>0$, then we have
\beeq\label{6-est-Bessel-Integral} \int^\infty_0 \frac{J_\mu(t)
J_\nu(t)}{t^\la} d t
=\frac{\Ga(\la)\Ga(\frac{\mu+\nu-\la+1}{2})}{2^\la
\Ga(\frac{\mu-\nu+\la+1}{2})
\Ga(\frac{\nu-\mu+\la+1}{2})\Ga(\frac{\mu+\nu+\la+1}{2})}.\eneq
\end{lem}

\section{Trace Lemma and Generalized Strichartz Estimates}\label{6-sec-traceStri}
In this section, we prove the trace lemma and the generalized
Morawetz estimates, in the presence of angular regularity. Moreover,
we give the proof of the generalized Strichartz estimates for the
wave equation (Theorem \ref{60-thm-Strichartz-Sterbenz}).

\subsection{Dual Verison of the Trace Lemma} 
\label{6-sec-trace} In this subsection, we prove the dual version of the trace lemma
stated in Theorem \ref{60-thm-trace}, which plays the central role
in this paper. For convenience, we restate it here.

\begin{thm}\label{6-thm-trace}
Let $b\in (1,n)$ and $n\ge 2$, then we have the following equivalent
relation
  \beeq\label{6-est-trace}
        \|\La_\omega^{\frac{b-1}{2}} |x|^{-\frac{b}{2}}\widehat{g d \sigma}(x)\|_{L^2_x}
        \simeq
        \|g\|_{L^2_\omega}.
  \eneq for any $g\in L^2_\omega$.
\end{thm}

By \eqref{6-est-trace}, we know that
  \beeq\label{6-est-trace2}
        \||x|^{-\frac{b}{2}}\widehat{g d \sigma}(x)\|_{L^2_x}
        \les
        \|g\|_{H^{s}_\omega}
  \eneq if  $s\ge \frac{1-b}{2}$ and $b\in(1,n)$.

\begin{rem}\label{6-rem-3.1}
For the estimate \eqref{6-est-trace} or \eqref{6-est-trace2}, the
condition on $b$ is also necessary. In fact, let $g=1$ and recalling
\eqref{6-est-F-d-sigma} and the well-known asymptotic of the Bessel
function \eqref{6-est-Bessel-asymptotic},
we know that we have
$$ \| |x|^{-\frac{b}{2}}\widehat{ d
 \sigma}(x)\|_{L^2_x}\simeq \| t^{-\frac{b+n-2}{2}} J_{\frac{n-2}{2}}(t)
 t^{{\frac{n-1}{2}}}\|_{L^2_{t>0}}
 \simeq  \| t^{-\frac{b-1}{2}} J_{\frac{n-2}{2}}(t) \|_{L^2_{t>0}} <\infty$$
 if and only if $b \in (1,n)$.
\end{rem}

Now we give the proof of Theorem \ref{6-thm-trace}.\\
\begin{prf}
For any
$g\in L^{2}_\omega$, we have the expansion formula with respect to
the spherical harmonics
$$g(\omega)=\sum_{k=0}^\infty \sum_{l=1}^{d(k)} a_{k,l} Y_{k,l}(\omega).$$
By \eqref{6-est-F-SphericalHarmonic2}, we have
$$\widehat{g d \sigma}(x)=
    \sum_k \sum_{l=1}^{d(k)} a_{k,l}
    (2 \pi)^{\frac{n}{2}} i^{- k}  r^{-\frac{n-2}{2}}
     J_{k+\frac{n-2}{2}}(r) Y_{k,l}(\omega)
:=\sum_k \sum_{l=1}^{d(k)} b_{k,l}(r) Y_{k,l}(\omega).$$ Since
$Y_{k,l}$ are ortho-normal bases on $L^2_\om$, we have
\begin{eqnarray*} \|\La_\omega^{\frac{b-1}{2}} |x|^{-\frac{b}{2}}
\widehat{g d \sigma}(x) \|_{L^2_x}^2 & = & \int^\infty_0  r^{-b}
\|\widehat{g d
\sigma}(r \omega)\|^2_{H^{\frac{b-1}{2}}_\omega} r^{n-1} d r \\
& \simeq & \sum_{k,l}
 \langle k\rangle^{b-1} \int^\infty_0  r^{-b}
|b_{k,l}(r)|^2 r^{n-1} d r \\
&=& \sum_{k,l} (2\pi)^n \langle k\rangle^{b-1} |a_{k,l}|^2
\|J_{k+\frac{n-2}{2}}(r) r^{\frac{1-b}{2}}\|_{L^2_r}^2.
\end{eqnarray*} If
$b\in(1,n)$, we can use the formula \eqref{6-est-Bessel-Integral}
with $\mu=\nu=k+\frac{n-2}{2}$ and $\la=b-1$ in the last integral.
Using the Stirling's formula \eqref{6-est-Gamma-asymptotic} for
large $t$, we have
$$\|J_{k+\frac{n-2}{2}}(r)
r^{\frac{1-b}{2}}\|_{L^2_r}^2=\frac{\Ga(b-1)\Ga(k+\frac{n-b}{2})}
{2^{b-1}\Ga(\frac{b}{2})^2 \Ga(k+\frac{n-2+b}{2})}\simeq \langle
k\rangle^{1-b}.$$ Thus, we get that
$$ \| \La_\omega^{\frac{b-1}{2}} |x|^{-\frac{b}{2}} \widehat{g d \sigma}(x)
\|_{L^2_x}^2 \simeq \sum_{k,l}  |a_{k,l}|^2 =
\|g\|^2_{L_{\omega}^{2}}.
$$\end{prf}

\subsection{Generalized Morawetz Estimates}\label{6-sec-Morawetz}

It is well-known that we can get the generalized Morawetz estimates
(or the local smoothing effect for the dispersive case $a >
1$) from the knowledge of the trace
lemma \eqref{6-est-trace2}. Here we can give a more refined estimate as stated in Theorem \ref{60-thm-Morawetz} because of the
improved trace lemma in Section \ref{6-sec-trace}.

{\noindent{\bf Proof of Theorem \ref{60-thm-Morawetz}.} }
Let $\xi=\la \om$ with $\om\in S^{n-1}$ and $s=\la^a$, we first
formally write $e^{i t \D^a} f$ as follows,
\begin{eqnarray*}
     e^{i t \D^a} f(x)
& =& (2 \pi)^{-n}\int e^{i (t |\xi|^a + x \cdot \xi)}     \hat{f}(\xi) d \xi\\
& =& (2 \pi)^{-n}\int_{S^{n-1}} \int_0^\infty  e^{i (t
     \la^a + \la x \cdot \om)}   \hat{f}(\la \om) \la^{n-1} d \la d \sigma(\om)\\
& =& \frac{1}{(2 \pi)^{n} a} \int_{S^{n-1}} \int_0^\infty  e^{i (t
     s + s^{\frac{1}{a}} x \cdot \om)}   \hat{f}(s^{\frac{1}{a}} \om) s^{\frac{n}{a}-1}
     d s d \sigma(\om)\\
& =& \frac{1}{(2 \pi)^{n} a} \int_0^\infty  e^{i t s}
     (\hat{f}(s^{\frac{1}{a}} \om) d \sigma (\omega))^{\vee}(s^{\frac{1}{a}} x )
     s^{\frac{n}{a}-1} d s\\
& =& \frac{1}{(2 \pi)^{n} a}
     \mathcal{F}^{-1}_s \{(\hat{f}(s^{\frac{1}{a}} \om) d \sigma (\omega))^{\vee}
     (s^{\frac{1}{a}} x )
     s^{\frac{n}{a}-1} H(s)\} (t) :=(h(\cdot, x))^\vee (t).
\end{eqnarray*}
By \eqref{6-est-trace} and the above expression, we know that
\begin{eqnarray*}
\||x|^{-\frac{b}{2}} e^{i t \D^a} f\|_{L^2_{t,x} } & \simeq &
\||x|^{-\frac{b}{2}} h(s, x)\|_{L^2_{s,x}} \\
& \simeq & \||x|^{-\frac{b}{2}} (\hat{f}(s^{\frac{1}{a}} \cdot) d
\sigma )^{\vee}( x ) s^{\frac{n+b}{2 a}-1} \|_{L^2_{s,x}}\\
& \simeq & \| \La_\omega^{\frac{1-b}{2}} \hat{f}(s^{\frac{1}{a}}
\om) s^{\frac{n+b-2 a}{2 a}} \|_{L^2_{s,\om}}\\
& \simeq & \| \La_\omega^{\frac{1-b}{2}} \hat{f}(s \om)
s^{\frac{n+b-2 a}{2}} s^{\frac{a-1}{2}} \|_{L^2_{s,\om}}\\
& \simeq & \| \La_\omega^{\frac{1-b}{2}} \hat{f}(\xi)
|\xi|^{\frac{b-a}{2}}  \|_{L^2_{\xi}} \simeq
\|\D^{\frac{b-a}{2}}\La_\omega^{\frac{1-b}{2}}f\|_{L^2_x}.
\end{eqnarray*} This is just the required estimate \eqref{60-est-Morawetz}.
Note that in the last equality, we have used the fact that
$$\|\La_\omega^{s} f\|_{L^2_x}\simeq \|\La_\omega^{s} \hat{f}\|_{L^2_x}.$$
This result can be easily seen,  write $f=\sum_{k,l}
f_{k,l}(r)Y_{k,l}(\omega)$, we have
$$\hat f(\xi)=\sum_{k,l} g_{k,l}(|\xi|)Y_{k,l}(\xi/|\xi|),$$
by \eqref{6-est-F-SphericalHarmonic}, and then applying
\eqref{6-est-Sobolev-Sphere},
\begin{eqnarray*}
  \|\La_\omega^{s} f\|_{L^2_x}&\simeq& \|\left<k\right>^s
f_{k,l}(|x|)\|_{l_{k,l}^2 L^2_x}\\
& \simeq& \|\left<k\right>^s f_{k,l}(|x|)Y_{k,l}(x/|x|)\|_{l_{k,l}^2
L^2_x}\\
& \simeq&  \|\left<k\right>^s
g_{k,l}(|\xi|)Y_{k,l}(\xi/|\xi|)\|_{l_{k,l}^2L^2_\xi}\\
& \simeq&  \|\left<k\right>^s
g_{k,l}(|\xi|)\|_{l_{k,l}^2L^2_\xi}\simeq \|\La_\omega^{s}
\hat{f}\|_{L^2_x}.
\end{eqnarray*}

If we apply \eqref{60-est-trace-v2} instead of \eqref{6-est-trace}
in the proof, we can get the estimate \eqref{60-est-Morawetz-Local}
by using a similar argument. More precisely, let $R_s=s^{1/a}R$, we
have
\begin{eqnarray*}
R^{-\frac{1}{2}} \| e^{i t \D^a} f\|_{L^2_{t,B(x_0,R)} } & \simeq &
R^{-\frac{1}{2}}\| h(s, x)\|_{L^2_{B(x_0,R), s}} \\
& \simeq & R_s^{-\frac{1}{2}} \| (\hat{f}(s^{\frac{1}{a}} \cdot) d
\sigma )^{\vee}( x ) s^{\frac{n+1}{2 a}-1} \|_{L^2_{s,B(x_0,R_s)}}\\
& \les & \| \hat{f}(s^{\frac{1}{a}}
\om) s^{\frac{n+1-2 a}{2 a}} \|_{L^2_{s,\om}}\\
& \simeq & \| \hat{f}(s \om)
s^{\frac{n+1-2 a}{2}} s^{\frac{a-1}{2}} \|_{L^2_{s,\om}}\\
& \simeq & \|  \hat{f}(\xi) |\xi|^{\frac{1-a}{2}}  \|_{L^2_{\xi}}
\simeq \|\D^{\frac{1-a}{2}} f\|_{L^2_x}.
\end{eqnarray*}
This completes the proof. {\hfill  {\vrule height6pt width6pt
depth0pt}\medskip}

\subsection{Generalized Strichartz Estimates for the Wave Equation}\label{6-sec-Strichartz}
In this subsection, we prove Theorem
\ref{60-thm-Strichartz-Sterbenz}, based on Proposition 3.4 and 3.5
in Sterbenz \cite{Stbz05}.

Let $(q_0,r_0)$ and $(q_1, r_1)$ be the endpoints of the classical
and generalized Strichartz estimates, i.e., $(q_1,r_1)=(2,
2\frac{n-1}{n-2})$ and
$$(q_0,r_0)=\left\{\begin{array}{ll}
(4,\infty)&n=2,\\
(2, 2\frac{n-1}{n-3})& n\ge 3.
\end{array}\right.$$

Set $q_\eta=2$ for $n\ge 3$ and $r_\eta=\infty$ for $n=2$.  We
recall first Proposition 3.4 and 3.5 in \cite{Stbz05}.
\begin{prop}\label{60-thm-sterbenzProp}
Let $u_{1,N}=e^{\pm i t \D} u_{1,N}(0)$ be a unit frequency, angular
frequency localized (at $2^N$) solution to the homogeneous wave
equation $\Box u_{1,N}=0$. Then for every $\eta>0$, there is a
$(q_\eta,r_\eta)$ such that $(q_\eta,r_\eta)\rightarrow(q_1,r_1)$ as
$\eta \rightarrow 0$ and the following estimate holds
\beeq\label{60-est-AngularSterProp} \|u_{1,N}\|_{L^{ q_\eta}_t
L^{r_\eta}_x} \le C_{\eta} N^{1/2+\eta} \|u_{1,N}(0)\|_{L^2_x}.
\eneq
\end{prop}

First, recall that to give the proof of Theorem
\ref{60-thm-Strichartz-Sterbenz}, we need only to prove the case
$r<\infty$, by the argument in \cite{FW2} (which uses generalized
Gargliardo-Nirenberg estimate and $L^q_t L^r_x$ estimate with
$r<\infty$ to control non-endpoint $L^q_t L_x^\infty$ norm, and can
essentially be viewed as a result of interpolation).

By the Littlewood-Paley theory and the Littlewood-Paley-Stein theory
(see Strichartz \cite{Stri72}), to prove Theorem
\ref{60-thm-Strichartz-Sterbenz} with $r<\infty$, we need only to
prove the following inequalities for $u_{1,N}$,
\beeq\label{60-est-AngularReduced} \|u_{1,N}\|_{L^{ q}_t L^r_x}\les
N^{s_{kn} + 2 \ep} \|u_{1,N}(0)\|_{L^2_x}\eneq for any $\ep>0$. By
interpolating with the trivial estimate $(q,r)=(\infty, 2)$, we can
further reduce it to the estimate with $r\in(r_1,r_0)$ and $q=2$ for
$n\ge 3$, or $q\in(q_1,q_0)$ and $r=\infty$ for $n=2$.

Recall that for the endpoint $(q_0,r_0)$, we have the following
estimate (see \cite{KeTa98} for $n\neq 3$ with $\ep=0$ and
\cite{MaNaNaOz05} for $n=3$) \beeq\label{60-est-AngularEndpt}
\|u_{1,N}\|_{L^{ q_0}_t L^{r_0}_x}\les N^{\ep}
\|u_{1,N}(0)\|_{L^2_x}\eneq for any $\ep>0$.

Now we can prove the required estimate by interpolation. Fix the
parameter $\ep>0$, we choose $\eta\ll 1$ to be fixed later, such
that $(q_\eta,r_\eta)\in [q_1, q]\times [r_1, r]$, and then choose
$t_\eta\in [0,1]$ such that
\beeq\label{60-est-AngularInterpo}\frac{t_\eta}{q_\eta}+\frac{1-t_\eta}{q_0}=\frac{1}{q},\
\frac{t_\eta}{r_\eta}+\frac{1-t_\eta}{r_0}=\frac{1}{r} \eneq i.e.,
$$t_\eta=(\frac{1}{q}-\frac{1}{q_0})/(\frac{1}{q_\eta}-\frac{1}{q_0})\textrm{ for }n=2 \textrm{ and }
t_\eta=(\frac{1}{r}-\frac{1}{r_0})/(\frac{1}{r_\eta}-\frac{1}{r_0})\textrm{
for }n\ge 3.
$$
Thus if we have \beeq\label{60-est-AngularCondition}
(\frac{1}{2}+\eta)t_\eta \le s_{kn} +\ep,\eneq then by
\eqref{60-est-AngularSterProp} and \eqref{60-est-AngularEndpt},
\begin{eqnarray*} \|u_{1,N}\|_{L^{ q}_t L^r_x}&\les&
\|u_{1,N}\|_{L^{
q_\eta}_t L^{r_\eta}_x}^{t_\eta}\|u_{1,N}\|_{L^{ q_0}_t L^{r_0}_x}^{1-t_\eta}\\
& \les& N^{(\frac{1}{2}+\eta)t_\eta+\ep (1-t_\eta)} \|u_{1,N}(0)\|_{L_x^2}\\
&\les& N^{s_{kn}+2 \ep} \|u_{1,N}(0)\|_{L_x^2},\end{eqnarray*} which
gives the required estimates \eqref{60-est-AngularReduced}.

Note that if $n=2$, we have
$$ \lim_{\eta\rightarrow 0} (\frac{1}{2}+\eta) t_\eta = \frac{1}{2}
(\frac{1}{q}-\frac{1}{q_0})/(\frac{1}{q_1}-\frac{1}{q_0})=\frac{2}{q}-\frac{1}{2}=s_{kn},$$
and if $n\ge 3$,
$$ \lim_{\eta\rightarrow 0} (\frac{1}{2}+\eta) t_\eta = \frac{1}{2}
(\frac{1}{r}-\frac{1}{r_0})/(\frac{1}{r_1}-\frac{1}{r_0})=(n-1)(\frac{1}{r}-\frac{1}{r_0})=s_{kn}.$$
Now we can choose $\eta$ sufficiently small, such that the estimate
\eqref{60-est-AngularCondition} holds true. This completes the proof
of Theorem \ref{60-thm-Strichartz-Sterbenz}.

\section{Strauss' Conjecture with rough data}
\label{6-sec-SLW}

Let $n\ge 2$, $F_p(u)=\la |u|^p$ ($\la\in \bR \backslash \{0\}$,
$p>1$), $s_c=\frac{n}{2}-\frac{2}{p-1}$, $p_{conf}=1+\frac{4}{n-1}$,
$p_h=1+\frac{4 n}{(n+1)(n-1)}$  and $p_c$ be the solution of the
quadratic equation
$$ (n-1) p_c^2 - (n+1) p_c - 2 = 0,\ p_c>1. $$

In this section, we apply the inequalities obtained in Section
\ref{6-sec-intro} and Section \ref{6-sec-traceStri} to the study of
the semi-linear wave equation \eqref{60-eqn-SLW} with small data
$(f,g)\in \dot{H}^{s_c,s_1}_\om \times \dot{H}^{s_c-1, s_1}_\om$ for
some $s_1$.

\subsection{Global results for $p_h<p<p_{conf}$}\label{6-sec-LdSo}
In this subsection, we give the proof of Theorem
\ref{60-thm-LdSo95}, following the arguments in Section 8 of
Lindblad-Sogge \cite{LdSo95}.

For the proof, we will invoke the generalized Strichartz estimates
in Theorem \ref{60-thm-Strichartz-Sterbenz} and the inhomogeneous
inequality of Harmse-Oberlin \eqref{60-est-Harmse}.

We use the usual contraction argument to give the proof. Let $q=
\frac{(n+1)(p-1)}{2}$ and
$$ X=\{u\in C_t \dot{H}_x^{s_c} \cap C^1_t \dot{H}_x^{s_c-1}\cap L_{t,x}^q : \|u\|_{L_{t,x}^q}\le C \ep\}$$
with $C>1$ to be determined, and $s_{kn}=\frac{1}{2}-s_c$. Define a
map $T: u \mapsto v$ for $u\in X$, such that $v$ is the solution of
the equation
$$(\pt^2-\Delta ) v = F_p (u), v(0,x)=f, \pt v(0,x)=g.$$
Hereafter, we denote by $v_{hom}$ and $v_{inh}$ the homogeneous and
inhomogeneous part of $v$ respectively, i.e., $v_{hom}=T(0)$ and
$v_{inh}=v-v_{hom}$.

Since $u\in X$, we have $F_p(u)\in L^{q/p}_{t,x}$. Note that
$p_h<p<p_{conf}$ if and only if
$\frac{n-1}{2(n+1)}<\frac{1}{q}<\frac{n-1}{2 n}$, which means that
the index pair $(q,q)$ satisfy \eqref{60-est-Conj1-wave}. Since
$\|(f,g)\|_{ \dot{H}^{s_c,s_1}_\om \times \dot{H}^{s_c-1,
s_1}_\om}\le \ep$ with $s_1>s_{kn}$, we get from Theorem
\ref{60-thm-Strichartz-Sterbenz}
$$\|v_{hom}\|_{
C_t \dot{H}^{s_c,s_1}_\om \cap L^q_{t,x}}+\|\pt v_{hom}\|_{C_t
\dot{H}^{s_c-1, s_1}_\om}\le C_1 \|(f,g)\|_{ \dot{H}^{s_c,s_1}_\om
\times \dot{H}^{s_c-1, s_1}_\om}\le  \frac{C}{2} \ep.$$

Note that $(q/p)'>2\frac{n+1}{n-1}$ if and only if $p<p_{conf}$, we
get that
$$\|v_{inh}\|_{C_t \dot{H}_x^{s_c} \cap C^1_t \dot{H}_x^{s_c-1} \cap L^q_{t,x}}
\le C_2 \|F_p(u)\|_{L^{q/p}_{t,x}} \le C_2 |\la| (C\ep)^p \le
\frac{C}{2}\ep$$ by the inhomogeneous version of the classical
Strichartz estimates \eqref{60-est-wave-Stri} and the inhomogeneous
estimate \eqref{60-est-Harmse} with $r=q/p$, if $\ep$ is
sufficiently small.

Similarly, one has $$\|v_1-v_2\|_{L^q_{t,x}}\le C_3
\|F_p(u_1)-F_p(u_2)\|_{L^{q/p}_{t,x}} \le C_4 (C\ep)^{p-1}
\|u_1-u_2\|_{L^q_{t,x}}\le \frac{1}{2} \|u_1-u_2\|_{L^q_{t,x}}.$$

Thus we have proved that for $\ep>0$ small enough, $T$ is a
contraction map on $X$, which completes the proof of Theorem
\ref{60-thm-LdSo95}.

\subsection{Strauss' Conjecture with mild rough data for $n\le
4$}\label{6-sec-strauss} In this subsection, we give the proof of
Theorem \ref{60-thm-Hidano} by using the weighted Strichartz
estimates \eqref{60-est-Morawetz-Strichartz} in Theorem
\ref{60-thm-Strichartz-weighted} and Sobolev inequality
\eqref{60-est-Sobolev-dual} in Corollary \ref{60-thm-Sobolev-trace}.

Since the method of proof is just the usual contraction argument (as
in subsection \ref{6-sec-LdSo}), we need only to give some of the
key inequalities here.

First, for $s_1=\frac{1}{p-1}$, $s_2=s_1+s_c-s_{sb}$,
$\al=\frac{n+1}{p}-\frac{2}{p-1}$, and any solution $u$ of the
equation $(\pt^2-\Delta) u=0$ with data $(f,g)\in
\dot{H}^{s_c,s_1}_\om \times \dot{H}^{s_c-1, s_1}_\om$, we have
$u\in C_t \dot{H}_{\om}^{s_c,s_1}\cap C^1_t
\dot{H}_{\om}^{s_c-1,s_1}$ and get from
\eqref{60-est-Morawetz-Strichartz} that
$$|x|^{-\al} u\in L^{p}_{t,|x|^{n-1}d|x|}
H^{s_2}_\om,$$ if $s_c-s_{sb}\in (0,\frac{n-1}{2})$, i.e., $p>p_c$.

Since $n\le 4$ and $p>p_c\ge 2$, we have $\frac{n-1}{2}< 2 \le [p]$,
where $[p]$ stand for the integer part of $p$. Note that by the
Moser estimate \eqref{6-est-Moser-Sphere} and the Sobolev embedding,
we get the following estimate on $S^{n-1}$ \beeq\label{6-est-Moser}
\|F_{p}(u)\|_{H^a_\om}\les \|u\|_{H^b_{\om}}^p \eneq with $b\ge
\frac{p-1}{p}\frac{n-1}{2}+\frac{a}{p}$ if $0\le a<\frac{n-1}{2}$
(and thus $a\le [p]$). Then by letting $b=s_2$ and $a=p
s_2-\frac{n-1}{2}(p-1)=\frac{n-1}{2}-\frac{1}{p-1} <\frac{n-1}{2}$,
we have
$$|x|^{-\al p} F_p(u)\in L^1_{t,|x|^{n-1}d|x|} H^{a}_\om.$$

Recall the estimate \eqref{60-est-Sobolev-dual}, if $\frac{1}{2}-s_c
\in (0, \frac{n-1}{2})$, i.e., $1+\frac{2}{n-1}<p<p_{conf}$, then
$$F_p(u)\in L^1_t \dot{H}^{s_c-1, a+\frac{1}{2}-s_c}_\om.$$ Note
that $a+\frac{1}{2}-s_c=\frac{1}{p-1}= s_1$, we have $F_p(u)\in
L^1_t \dot{H}^{s_c-1, s_1}_\om$.

By the classical energy estimate, we can get again $(v, \pt v)\in
C_t (\dot{H}_{\om}^{s_c,s_1}\times \dot{H}_{\om}^{s_c-1,s_1})$
 and
$$|x|^{-\al} v\in L^{p}_{t,|x|^{n-1}d|x|}
H^{s_2}_\om,$$ if $(\pt^2-\Delta) v=F_p(u)$ with initial data
$(f,g)$.

 If we define the solution space $X$ to be
$$ X=\{u\in C_t \dot{H}_{\om}^{s_c,s_1}\cap C^1_t \dot{H}_{\om}^{s_c-1,s_1} :
\||x|^{-\al} u\|_{ L^{p}_{t,|x|^{n-1}d|x|} H^{s_2}_\om}\le C
\ep\},$$ then, at last, combining all above, we can prove that the
problem \eqref{60-eqn-SLW} with $p_c<p<p_{conf}$ is global
well-posed for small data in the space $C_t
\dot{H}_{\om}^{s_c,s_1}\cap C^1_t \dot{H}_{\om}^{s_c-1,s_1}$ with
$s_1= \frac{1}{p-1}$. Thus we get Theorem \ref{60-thm-Hidano}.

\section{An Application to the Schr\"{o}dinger Equation}\label{6-sec-SLS}
In this section, we prove Theorem \ref{60-thm-SLS}. As in Section
\ref{6-sec-strauss}, we will only give the key inequalities for the
proof.

Let $n\ge 3$, $1<p< p_{L2}$, $s_1=\frac{1}{p-1}$,
$\al=\frac{n+2}{q}-\frac{2}{p-1}$ and $q\ge 2$ to be determined
later. The solution space is $$ X=\{u\in C_t \dot{H}_{\om}^{s_c,s_1}
: 
\||x|^{-\al} u\|_{
L^{q}_{t,|x|^{n-1}d|x|} H^{s_2}_\om}\le C \ep\},$$

By \eqref{60-est-Morawetz-Strichartz2},  we have $e^{-i t \Delta}
f\in X$ with \beeq\label{60-est-SLS-CondS1}
s_2=s_1+\frac{n-1}{2}+\al-\frac{n}{q}=\frac{n-1}{2}+\frac{2}{q}-\frac{1}{p-1},\eneq
for any $f \in \dot{H}^{s_c,s_1}_\om$ with norm $\le \ep$, if
$\frac{n}{q}-\al=\frac{2}{p-1}-\frac{2}{q}\in (0, \frac{n-1}{2})$,
i.e., \beeq\label{60-est-SLS-CondB1}
\frac{2}{p-1}-\frac{n-1}{2}<\frac{2}{q}<\frac{2}{p-1}.\eneq

By the Moser estimate \eqref{6-est-Moser}, if $s_3=p
s_2-\frac{n-1}{2}(p-1)\in [0,[p]]$ and $s_3<\frac{n-1}{2}$, we have
\beeq |x|^{-\al p} F_p(u)\in L^{{q/p}}_{t,|x|^{n-1}d|x|} H^{s_3}_\om
\eneq for any $u$ such that $|x|^{-\al} u\in L^{q}_{t,|x|^{n-1}d|x|}
H^{s_2}_\om$. The restrictions on $s_2$ and $s_3$ can be  reduced to
that on $q$, i.e.
\beeq\label{60-est-sls-condq}\frac{1}{p-1}-\frac{n-1}{2 p}\le
\frac{2}{q} < \frac{1}{p-1}-\frac{n-3}{2 p}.\eneq

For any $u\in X$, define $v=T u$ to  be the solution of the equation
$(i\pt-\Delta) v=F_p(u)$ with $v(0)=0$.  We want to show that $v\in
X$. Recall the estimate \eqref{60-est-Mora-Stri-dual2}, we have
$v\in C_t \dot{H}_{\om}^{s_c,s_1}$, if $q/p\le 2$, $s_1\le
s_3+\frac{n-1}{2}-\al p-\frac{n}{(q/p)'}$ and $\frac{n}{(q/p)'}+\al
p\in (0, \frac{n-1}{2})$, i.e., \beeq\label{60-est-sls-condS3}
s_1\le s_3-s_c+\frac{3}{2}-\frac{2p}{q},\eneq
\beeq\label{60-est-SLS-condb2} \frac{2}{q}\ge \frac{1}{p}\textrm{,
and }
\frac{2}{p-1}-\frac{n}{p}<\frac{2}{q}<\frac{2}{p-1}-\frac{n+1}{2p}.\eneq
Moreover, by estimate \eqref{60-est-Mora-Stri-Inhom}, we have
$$|x|^{-\al} v\in L^{q}_{t,|x|^{n-1}d|x|} H^{s_2}_\om$$ if
$s_2-s_3\le (\frac{n-1}{2}+\al-\frac{n}{q})+(\frac{n-1}{2}-\al
p-\frac{n}{(q/p)'})$, i.e., \beeq\label{60-est-sls-condS4}
s_2-s_3\le
(\frac{2}{q}-\frac{1}{2})+(\frac{2}{(q/p)'}-\frac{1}{2})=1-2\frac{p-1}{q}.\eneq

The above information is sufficient for us to prove the Theorem
\ref{60-thm-SLS} by using the contraction argument. Thus to complete
the proof, we need only to check the requirement on $q$. By the
expression of $s_2$ and $s_3$, it is easy to check that the
inequality \eqref{60-est-sls-condS3} and \eqref{60-est-sls-condS4}
hold with equality.

The restrictions on $q\ge 2$ is just \eqref{60-est-SLS-CondB1},
\eqref{60-est-sls-condq} and \eqref{60-est-SLS-condb2}. Recall that
we assume $p<p_{L2}<\frac{2n}{n-1}$, then
$$\frac{2}{p-1}-\frac{n}{ p}<\frac{2}{p-1}-\frac{n-1}{2}.$$ This
means that the restrictions on $q$ is just
\beeq\label{60-est-sls-resu-q0} \frac{2}{q} \in [\frac{1}{p}, 1]\cap
(\frac{2}{p-1}-\frac{n-1}{2}, \frac{2}{p-1}-\frac{n+1}{2 p})\cap
[\frac{1}{p-1}-\frac{n-1}{2 p}, \frac{1}{p-1}-\frac{n-3}{2 p}).\eneq

This is our requirement on $q$. To complete the proof of Theorem
\ref{60-thm-SLS}, we conclude that this set is nonempty if \beeq
\label{60-est-sls-resu-p} p-1\in (\sqrt{\frac{2}{n-1}},
\frac{4}{n})\textrm{ and }n \le 6,\eneq i.e. $p_l<p<p_{L2}$ and $n
\le 6$. This can be directly calculated. The somewhat strange index
$p_l$ comes from the condition
$$\frac{2}{p-1}-\frac{n-1}{2}<\frac{1}{p-1}-\frac{n-3}{2 p}.$$
We want $$\frac{1}{p}<\frac{1}{p-1}-\frac{n-3}{2 p}\Leftrightarrow
p<1+\frac{2}{n-3}$$ is true for any $p=p_{L2}-\ep$ with $\ep$ small
enough, which gives us the restriction on the dimension $n\le 6$.
This finishes the proof of Theorem \ref{60-thm-SLS}.

\medskip
\noindent{\sl Acknowledgement.} We would like to thank the referee
for a careful and thorough reading of the manuscript and a number of
criticisms and helpful suggestions.

\end{document}